\documentclass[a4paper,twoside]{article}
\usepackage{amsmath,amssymb,amsfonts,amsthm}
\pdfoutput=1
\usepackage{hyperref}
\usepackage{diagrams}
\begin{document}
\title{A review of some recent work on hypercyclicity\footnote{Invited paper,
  Workshop celebrating the 65$^{th}$ of L. A. Cordero, Santiago de Compostela, June 27-29, 2012}}
\author{  C.T.J. Dodson\\
School of  Mathematics, University of Manchester\\
Manchester M13 9PL UK.\\
{\small ctdodson@manchester.ac.uk}}
\pagestyle{myheadings}
\markboth{Review of recent work on hypercyclicity}{C.T.J. Dodson}
\date{}
\maketitle\maketitle\addcontentsline{toc}{section}{Title and Abstract}

%===================================================================
%\theoremstyle{definition}
\newtheorem{definition}{Definition}[section]
\newtheorem{theorem}{Theorem}[section]
\newtheorem{proposition}[theorem]{Proposition}
\newtheorem{corollary}[theorem]{Corollary}
\newtheorem{lemma}[theorem]{Lemma}
\newtheorem{remark}[theorem]{Remark}
%===================================================================
\begin{abstract}
 \noindent Even linear operators on infinite-dimensional spaces can display interesting dynamical properties and yield important links among functional analysis, differential and global geometry and dynamical systems, with a wide range of applications. In particular, hypercyclicity is an essentially infinite-dimensional property, when iterations of the operator generate a dense subspace. A Fr\'{e}chet space admits a hypercyclic operator if and only if it
is separable and infinite-dimensional. However, by considering the semigroups generated by multiples of operators, it is possible to obtain hypercyclic behaviour on finite dimensional spaces. This article gives a brief review of some recent work on hypercyclicity of operators on Banach, Hilbert and Fr\'{e}chet spaces.

\noindent{\bf MSC:} 58B25 58A05 47A16, 47B37\\
{\bf Keywords:}  Banach space, Hilbert space, Fr\'{e}chet space, bundles,  hypercyclicity.
\end{abstract}
\section{Introduction}
In a number of cases that have significance in global analysis~\cite{Hamilton,Neeb}, physical field theory~\cite{Smolentsev}, dynamical systems~\cite{BayartMath,ShapiroNotes01,GrE11} and  finance theory~\cite{EmamiradGG11}, Banach space representations may break down and we need Fr\'{e}chet spaces, which have weaker requirements for their topology.  Fr\'{e}chet spaces of sections arise naturally as configurations of a physical field where the moduli space, consisting of inequivalent configurations of the physical field, is the quotient of the infinite-dimensional configuration space $\mathcal{X}$ by the appropriate
symmetry gauge group. Typically, $\mathcal{X}$ is modelled on a Fr\'{e}chet
space of smooth sections of a vector bundle over a closed manifold. Countable products
of an infinite-dimensional Banach space are non-normable Fr\'{e}chet spaces. See the notes of
Doma\'{n}ski~\cite{Domanski10} for a collection of results on spaces of analytic functions and
linear operators on them, including projective limit spaces.

Smolentsev~\cite{Smolentsev} and Clarke~\cite{Clarke} discuss the metric geometry of the Fr\'{e}chet manifold of all $C^\infty$ Riemannian metrics on a fixed closed finite-dimensional orientable manifold. Micheli et al~\cite{MMM12} discuss Sobolev metrics and geodesic behaviour on groups of diffeomorphisms of a finite-dimensional manifold under the condition that the diffeomorphisms decay suitably rapidly to the identity.
Omori~\cite{Omori1,Omori4} provides further discussion of Lie-Fr\'{e}chet groups of diffeomorphisms of closed Riemannian manifolds as ILH-manifolds, that is as inverse (or projective) limits of Hilbert manifolds; unlike Fr\'{e}chet manifolds, Hilbert manifolds do support the main theorems of calculus.

These weaker structural constraints raise other problems: Fr\'{e}chet spaces lack a general solvability theory of differential equations, even linear ones; also, the space of continuous linear mappings does not remain in the category while the space of linear isomorphisms does not admit a reasonable Lie group structure. Such shortcomings can be worked round to a certain extent by representing Fr\'{e}chet spaces as projective limits of Banach spaces and in the manifold cases by requiring the geometric structures to carry  through as projective limits, see Galanis et al.~\cite{Gal2,VG 1,Dod-Gal1,DGV 1,DGV 2} for results on tangent and second tangent bundles, frame bundles and generalized Lie groups, cf.~\cite{Dodson11} for a survey. In a detailed study of Lie group actions on Banach spaces, with several appendices on the necessary calculus, Walter~\cite{Walter10} elaborated details of solutions
of differential equations on each step of a projective limit and integration of some Lie algebras of vector fields.

An open problem is the extension to Banach, Hilbert and Fr\'{e}chet bundles of the results on projection and lifting of harmonicity for tangent, second tangent and frame bundles obtained with Vazquez-Abal~\cite{DVBolIt,DVMatZem}\footnote{Which work originally brought the author to Santiago at the invitation of Luis Cordero in the early 1980s}, for finite-dimensional Riemannian manifolds:
\begin{diagram}
   (FM,Fg) & \rTo^{\pi_{FM}} & (M,g)& \lTo^{\pi_{TM}} & (TM,Tg)\\
  \dTo^{Ff}    &    &  \dTo_{f}  &    &  \dTo_{Tf}\\
   (FN,Fh)               &  \rTo^{\pi_{FN}}    &      (N,h)   &  \lTo^{\pi_{TN}}    & (TN,Th)
\end{diagram}
In this diagram $f$ needs to be a local diffeomorphism of Riemannian manifolds for the frame bundle morphism $Ff$ to be defined. It was shown that $Ff$ is totally geodesic if and only if $f$ is totally geodesic; when $f$ is a local diffeomorphism of flat manifolds then $Ff$ is harmonic if $f$ is harmonic. Also, the diagonal map $\pi_{FN}\circ Ff=f\circ \pi_{FM}$ is harmonic if and only if $f$ is harmonic, and $Ff$ is harmonic if and only if $Tf$ is harmonic. Sanini~\cite{Sanini} had already established the corresponding result for the tangent bundle projection: $Tf$ is totally geodesic if and only if $f$ is totally geodesic. It follows~\cite{DVMatZem}, using Smith~\cite{Smith}, that $\pi_{TM}$ is a harmonic Riemannian submersion and the diagonal map $\pi_{TN}\circ Tf=f\circ \pi_{TM}$ is harmonic if and only if $f$ is harmonic.

It would, for example, be interesting to extend the above to the infinite dimensional case of an inverse limit Hilbert (ILH) manifold $\mathbb{E} = \lim_{\infty\leftarrow s}\mathbb{E}^s,$
 of a projective system of smooth Hilbert manifolds $\mathbb{E}^s$, consisting of sections of a tensor bundle over a smooth compact finite dimensional Riemannian manifold $(M,g).$ Such spaces arise in geometry and physical field theory and they have many desirable properties but it is necessary to establish existence of the projective limits for various geometric objects. Smolentsev~\cite{Smolentsev} gives a detailed account of the underlying theory we
 need---that paper is particularly concerned with the manifold of sections of the bundle of smooth symmetric 2-forms on $M$ and its critical points for important geometric functionals.
 We may mention the work of Bellomonte and Trapani~\cite{BellomonteT12}who investigated directed systems of Hilbert spaces whose extreme spaces are the projective and the inductive limit of a directed contractive family of Hilbert spaces.

 Via the volume form on ($n$-dimensional compact) $(M,g)$ a weak induced metric on the space of tensor fields is $\int_M g(X,Y)$ but there is a stronger family~\cite{Smolentsev} of inner products on $\mathbb{E}^s,$ the completion Hilbert space of sections. For sections $X,Y$ of the given tensor bundle over $M$ we put
\begin{equation}\label{ip}
    (X,Y)_{g,s} = \sum_{i=0}^s\int_M g(\nabla^{(i)}X,\nabla^{(i)}Y) \ \ \ s\geq 0.
\end{equation}
Then the limit
 $\mathbb{E} = \lim_{\infty\leftarrow s}\mathbb{E}^s$
with limiting inner product $g_\mathbb{E}$ is a Fr\'{e}chet space with topology independent of the choice of metric $g$ on $M.$ In particular it is known, for example see  Omori~\cite{Omori1,Omori4} and Smolentsev~\cite{Smolentsev}, that the  smooth diffeomorphisms $f:(M,g)\rightarrow (M,g)$ form a strong ILH-Lie group ${\cal D}iff_M$ modelled on the ILH manifold
$$\Gamma(TM) = \lim_{\infty\leftarrow s} \Gamma^s(TM)$$
 of smooth sections of the tangent bundle. Moreover, the curvature and Ricci tensors are equivariant under the action of ${\cal D}iff_M,$ which yields the Bianchi identities as consequences. The diagram of Hilbert manifolds of sections of vector bundles over smooth compact finite dimensional Riemannian manifolds $(M,g),\  (N,h)$ with $\mathbb{E}=\Gamma(TM),\ \mathbb{F}=\Gamma(TN).$ Diagonal lift metrics are induced via the horizontal-vertical splittings defined by the Levi-Civita connections $\nabla^g, \nabla^h$ on the base manifolds (cf.~\cite{Sasaki,Domb62,Kow71}), effectively applying the required evaluation to corresponding projections; we abbreviate these to  $Tg_\mathbb{E}=(g_\mathbb{E},g_\mathbb{E}),\   Th_\mathbb{F}=(h_\mathbb{F},h_\mathbb{F}),$
\begin{diagram}
   (T\mathbb{E},Tg_\mathbb{E}) & \rTo^{T\phi} & (T\mathbb{F},Th_\mathbb{F})\\
  \dTo^{\pi_{T\mathbb{E}}}   &  &     \dTo^{\pi_{T\mathbb{F}}}   \\
    (\mathbb{E},g_\mathbb{E})   & \rTo^{\phi}    &      (\mathbb{F},h_\mathbb{F})
\end{diagram}
For example, a smooth map of Riemannian manifolds $f:(M,g)\rightarrow (N,h)$
defines a fibre preserving map $f^*$ between their tensor bundles and induces such a smooth map $\phi$ between the spaces of sections.
The Laplacian $\triangle$ on our Hilbert manifold $\mathbb{E}$ is defined by $\triangle=-{\rm div}\nabla^{\mathbb{E}} d$ where the generalized divergence $-{\rm div}$ is the trace of the covariant derivation operator $\nabla^{\mathbb{E}},$
so ${\rm div}$ is the adjoint of the covariant derivation operator $\nabla^{\mathbb{E}}.$ At this juncture we defer to  future studies the investigation of lifting and projection of harmonicity in ILH manifolds and turn to the characterization of linear operators
then review work reported in the last few years on the particular property of hypercyclicity, when iterations generate dense subsets.

\section{Dynamics of linear operator equations}
A common problem in applications of linear models is the characterization and solution of continuous linear operator equations on Hilbert, Banach and Fr\'{e}chet spaces. However, there are many open problems. For example, it is known that for a continuous linear operator $T$ on a Banach space $\mathbb{E}$ there is no non-trivial closed subspace nor non-trivial closed subset $A\subset\mathbb{E}$ with $TA\subset A,$ but this is an unsolved problem on Hilbert and Fr\'{e}chet spaces, cf. Martin~\cite{Martin11} and Banos~\cite{Banos11} for more discussion of invariant subspace problems. Shapiro's notes~\cite{ShapiroNotes01}  illustrate how continuous linear transformations of infinite dimensional topological vector spaces can have interesting dynamical properties, with new links among the theories of dynamical
systems, linear operators, and analytic functions. The notes of Doma\'{n}ski~\cite{Domanski10}
collect a wide range of results on commonly studied spaces of real analytic functions and linear operators on them. Also, his paper~\cite{Domanski10a} on the real analytic parameter dependence of solutions of linear partial differential equations has detailed solutions for a wide range of equations, establishing also a characterization of surjectivity of tensor products of general surjective linear operators on a wide class of spaces containing most of the natural spaces of classical analysis.

There has been
substantial interest from differential geometry and dynamical systems in hypercyclic operators, whose iterations generate dense subsets.
In this survey we look at some of the results on hypercyclicity of operators that have been reported in the last few years.

\section{Hypercyclicity properties}
A continuous linear operator  $T$ on a topological vector space  $\mathbb{E}$ is {\em hypercyclic} if, for some $f\in\mathbb{E}$, called a {\em hypercyclic vector}, the set $\{T^nf, n\geq 0 \}$ is dense in $\mathbb{E},$  and {\em supercyclic} if the projective space orbit
$\{ \lambda T^nf, \lambda \in \mathbb{C}, n\geq 0 \}$ is dense in $\mathbb{E}.$ These properties are called {\em weakly hypercyclic}, {\em weakly supercyclic} respectively, if $T$  has the property with respect to the weak topology---the smallest
topology for the space such that every member of the dual space is continuous
with respect to that topology. See the earlier reviews by Grosse-Erdmann~\cite{GrE,GrE03} and the recent books by Grosse-Erdmann and Manguillot~\cite{GrE11} and Bayart and Matheron~\cite{BayartMath} for more details of the development of the theory of hypercyclic operators.

If $T$ is invertible, then it is hypercyclic if and only if $T^{-1}$ is hypercyclic.
It is known for $\ell^p(\mathbb{N}),$ the Banach space of complex sequences with $p$-summable modulus $p\geq 1$ and backward shift operator
$B_{-1}:(x_0,x_1,x_2,\ldots)\mapsto (x_1,x_2,x_3,\ldots),$ that
 $\lambda B_{-1}$ is hypercyclic
on $\ell^p(\mathbb{N})$ if and only if $|\lambda|> 1.$
De La Rosa~\cite{delaRosa11}discussed operators which
are weakly hypercyclic, summarizing properties shared with hypercyclic operators, and proved the following about a weak hypercyclic $T:$\\
(i) $T\oplus T$ need not be weakly hypercyclic, with an example on , $\ell^p(\mathbb{N})\oplus\ell^p(\mathbb{N}), \ 1\leq p <\infty$\\
(ii)  $T^n$ is weakly hypercyclic for every $n>1$ \\
(iii) For all unimodular $\lambda\in \mathbb{C},$ we have $\lambda T$ weakly hypercyclic.\\
Thus, a weakly hypercyclic operator has many of the same properties
as a hypercyclic operator. For example, its adjoint, has no eigenvalue and every component of its spectrum must intersect the unit circle. However, De La Rosa~\cite{delaRosa11} \S3 summarized some known examples illustrating differences. Clements~\cite{Clements12} analyzed in detail the spectrum
for hypercyclic operators on a Banach space. Shkarin~\cite{Shkarin12} established a new 
criterion of weak hypercyclicity of a bounded linear operator on a Banach space. 

It was known from Godefroy and Shapiro~\cite{GodefroyShapiro91} that on every separable Banach space, hypercyclicity is equivalent to transitivity: i.e.  for every pair
of nonempty, norm open sets $(U,V),$ we have $T^n(U) \bigcap V\neq\emptyset$ for some $n\in \mathbb{N},$  and in particular, on the Fr\'{e}chet space of analytic functions on $\mathbb{C}^N$ every linear partial differential operator with constant coefficients and positive order has a hypercyclic vector. However, that proof does not carry over to the weak topology. Chan~\cite{Chan02} showed that on a separable infinite-dimensional complex Hilbert space $\mathbb{H}$ the set of hypercyclic operators is dense in the strong operator topology, and moreover the linear span of hypercyclic operators is dense in
the operator norm topology. The non-hypercyclic operators are dense in the set of bounded operators $B(\mathbb{H})$ on $\mathbb{H}$, but the hypercyclic
operators are not dense in the complement of the closed unit ball of $B(\mathbb{H})~\cite{Chan02}.$
Rezai~\cite{Rezaei11} investigated transitivity of  linear operators acting on a reflexive Banach space $\mathbb{E}$ with the weak topology. It was shown that a bounded operator, transitive on an open bounded subset of $\mathbb{E}$ with the weak topology, is weakly hypercyclic.

Evidently, if a linear operator is hypercyclic, then having a hypercyclic vector means also that it possesses a dense subspace in which all nonzero vectors are hypercyclic. A {\em hypercyclic subspace} for a linear operator is an infinite-dimensional closed subspace all of whose nonzero vectors are hypercyclic. Menet~\cite{Menet12} gave a simple critrion for a Fr\'{e}chet space 
with a continuous norm to have no hypercyclic subspaces; also if $P$ is a non-constant polynomial 
and $D$ is differentiation on the space of entire functions then $P(D)$ possesses a hypercyclic subspace.

On the Fr\'{e}chet space $\mathbb{H}(\mathbb{C})$ of functions analytic on $\mathbb{C},$ the translation by a fixed nonzero $\alpha\in \mathbb{C}$ is hypercyclic and so is the differentiation operator $f\mapsto f'.$ Ansari~\cite{Ansari97} proved that all infinite-dimensional separable Banach spaces admit hypercyclic operators. On the other hand, no finite-dimensional
Banach space admits a hypercyclic operator. Every nonzero power $T^m$ of a hypercyclic linear operator $T$ is hypercyclic, Ansari~\cite{Ansari}. Salas~\cite{Salas95} used backward weighted shifts on $\ell^2$ such that $T(e_i) = w_ie_{i-1} \ (i\geq 1) \ {\rm and} \  T(e_0)=0$ with positive $w_i$ to show that $T+I$ is hypercyclic.
 All infinite-dimensional separable Banach spaces admit
hypercyclic operators by Ansari~\cite{Ansari97}; however, Kitai~\cite{Kitai} showed that finite dimensional spaces do not. In particular a Fr\'{e}chet space admits a hypercyclic operator if and only if it is separable and infinite-dimensional and the spectrum of a hypercyclic operator must meet the unit circle. Bakkali and Tajmouati~\cite{BakkaliT12} have provided some further Weyl and Browder spectral characterizations of hypercyclic and supercyclic operators on separable Banach and Hilbert spaces.

A {\em sequence} of linear operators $\{T_n\}$ on $\mathbb{E}$ is called hypercyclic if, for some $f\in \mathbb{E},$ the set $\{T_nf, n\in \mathbb{N}\}$ is dense in $\mathbb{E};$ see Chen and Shaw~\cite{ChenShaw} for a discussion of related properties. The sequence $\{T_n\}$ is said to satisfy the {\em Hypercyclicity Criterion} for an increasing sequence $\{n(k)\}\subset \mathbb{N}$ if there are dense subsets $X_0,Y_0\subset \mathbb{E}$ satisfying (cf. also Godefroy-Shapiro~\cite{GodefroyShapiro91}):\\
{\bf Hypercyclicity  Criterion}
\begin{eqnarray}\label{HC}
   &(\forall f\in X_0)& T_{n(k)} f\rightarrow 0 \nonumber \\
   &(\forall g\in Y_0)& {\rm there \ is \ a \ sequence} \ \{u(k)\}\subset \mathbb{E} \
    {\rm such \ that} \ u(k) \rightarrow 0 \ {\rm and} \ T_{n(k)} u(k)\rightarrow 0. \nonumber
\end{eqnarray}
Bes and Peris~\cite{BesPeris} proved that on a  separable Fr\'{e}chet space $\mathbb{F}$ a continuous linear operator $T$ satisfies the Hypercyclicity Criterion if and only if $T\oplus T$ is hypercyclic on $\mathbb{F}\oplus\mathbb{F}.$  Moreover, if $T$ satisfies the Hypercyclicity Criterion then so does every power $T^n$ for $n\in \mathbb{N}.$ Rezai~\cite{Rezaei12} showed that 
such a $T$ with respect to a syndetic sequence (increasing positive integers $n_k$ with bounded 
$\sup(n_{k+1}-n_k)$) then $T$ satisfies the Kitai Criterion~\cite{Kitai}.  

A vector $x$ is called {\em universal} for a sequence of operators $(T_n: n\in \mathbb{N}\}$ on a Banach space $\mathbb{E}$ if $\{T_nx : n\in \mathbb{N}\}$ is dense; $x$ is called {\em frequently universal} if for each
non-empty open set $U\subset \mathbb{E}$ the set $K=\{n : T_n\in U\}$ has positive lower density,
namely $\liminf_{N\rightarrow\infty}|\{n\leq N: n\in K\} |/N >0.$
A {\em frequently hypercyclic} vector of $T$ is such that, for each
non-empty open set $U,$ the set $\{n : T^n\in U\}$ has positive lower density, a stronger
requirement than hypercyclicity. Drasin and Saksman~\cite{DrasinSaksman12} deduce optimal growth 
properties of entire functions frequently hypercyclic on the differentiation operator, cf. also 
Blasco et al.~\cite{BlascoEtAl10}. 
Bonilla and Grosse-Erdmann~\cite{BonillaGE07} extended a sufficient condition for frequent hypercyclicity from Bayart and Grivaux~\cite{BayartGrivaux06}, to frequent universality.
Beise~\cite{Beise12} extended this work and gave a sufficient condition for frequent universality in the Fr\'{e}chet case.

Extending the work of Yousefi and Rezaei~\cite{YR07},
Chen and Zhou~\cite{ChenZhou11} obtained necessary and sufficient conditions for the hypercyclicity of weighted composition operators (cf. also Bonet and Doma\'{n}ski~\cite{BonetDomanski11}) acting on the complete vector space of
holomorphic functions on the open unit ball $B_{N}$ of $\mathbb{C}^{N}.$ The weighted composition operators are constructed as follows.
Let $\varphi$ be a holomorphic self-map of $B_{N}$ then the composition operator with symbol $\varphi$ is $C_\varphi: f\mapsto f\circ\varphi$ for $f\in H(B_{N})$ the space of holomorphic maps on $B_{N}.$ The multiplication operator induced by $\psi\in H(B_{N})$ is $M_\psi(f)=\psi\cdot f$ and the weighted composition operator induced by $\psi,\varphi$ is $W_{\psi,\varphi}=M_\psi C_\varphi.$ Further results established that if $C_\varphi$ is hypercyclic then so is $\lambda C_\varphi$ for all unimodular $\lambda\in \mathbb{C};$ also, if $\varphi$ has an interior fixed point $w$ and $\psi\in H(B_{N})$ satisfies
$$|\psi(w)| < 1 < \lim_{|z|\rightarrow 1} \inf|\psi(z)|,$$
then the adjoint $W^*_{\psi,\varphi}$ is hypercyclic.

Zaj\c{a}c~\cite{Zajac12} characterized hypercyclic composition operators $C_\varphi:f\mapsto f\circ\varphi$ on the space of functions holomorphic on $\Omega\subset\mathbb{C}^N,$  a pseudoconvex domain and $\varphi$ is a holomorphic self-mapping of $\Omega.$ In the case when all the balls with respect to the Carath\'{e}odory pseudodistance are relatively compact in $\Omega$, he showed that much simpler characterization is possible (e.g. strictly pseudoconvex domains, bounded convex domains). Also, in such a class of domains, and in simply connected or infinitely connected planar domains, hypercyclicity of $C_\varphi$ implies it is hereditarily hypercyclic, i.e.  $C_\varphi\oplus C_\varphi$ is hypercyclic~\cite{BesPeris}.

Montes-Rodriguez et al.~\cite{MRRMS} studied the Volterra composition operators $V_\varphi$ for $\varphi$ a measurable self-map of $[0,1]$ on functions $f\in L^p[0,1], \ 1\leq p \leq \infty$
\begin{equation}\label{Volt}
    (V_\varphi f)(x) =\int_0^\varphi(x) f(t) dt
\end{equation}
These operators generalize the classical Volterra operator $V$ which is the case when $\varphi$ is the identity. $V_\varphi$ is measurable, and compact on $L^p[0,1].$
Consider the Fr\'{e}chet space  $\mathbb{F}=C_0[0,1),$ of continuous functions vanishing at zero with the topology of uniform convergence on compact subsets of $[0,1).$ It was known that the action of $V_\varphi$ on $C_0[0,1)$ is hypercyclic when $\varphi(x)=x^b, b\in (0,1)$ by Herzog and Weber~\cite{Herzog}. This result has now been extended by Montes-Rodriguez et al. to give the following complete characterization.
 \begin{theorem}~\cite{MRRMS}
For $\varphi\in C_0[0,1)$ the following are equivalent\\
{\bf (i)} $\varphi$ is strictly increasing with $\varphi(x)>x$ for $x\in (0,1)$\\
{\bf (ii)}  $V_\varphi$ is weakly hypercyclic\\
{\bf (iii)} $V_\varphi$ is hypercyclic.
\end{theorem}
Extending the work of Salas~\cite{Salas95,Salas99}, Montes-Rodriguez et al.~\cite{MRRMS} proved for every strictly increasing $\varphi$ with  $\varphi(x)< x, \ x\in(0,1]$
that $V_\varphi$ is supercyclic and $I+V_\varphi$ is hypercyclic when $V_\varphi$ acts on $L^p[0,1], \ p\geq 1, \ {\rm or \ on} \  \cal{C}_0[0,1].$  Shu et al.~\cite{ShuZhaoZhou11}  showed that the conjugate set $\{L^{-1}T L: L \ {\rm invertible}\}$ of any supercyclic operator
$T$ on a separable, infinite dimensional Banach space contains a path of supercyclic
operators which is dense with the strong operator topology, and the set of common
supercyclic vectors for the path is a dense $G_\delta$ set (countable intersection of open and
dense sets) if $\sigma_p(T^*)$ is empty.

Karami et al~\cite{Karami} gave examples of hypercyclic operators on $H_{bc}(\mathbb{E}),$ the space of bounded functions on compact subsets of Banach space $\mathbb{E}.$ For example, when  $\mathbb{E}$  has separable dual $\mathbb{E}^*$ then for nonzero $\alpha\in\mathbb{E},$ $T_\alpha:f(x)\mapsto f(x+\alpha)$ is hypercyclic.
As for other cases of hypercyclic operators on Banach spaces, it would be interesting to know when the property persists to projective limits of the domain space.

Yousefi and Ahmadian~\cite{YA} studied the case that $T$ is a continuous linear operator on an infinite dimensional Hilbert space $\mathbb{H}$ and left multiplication is hypercyclic with respect to the strong operator topology. Then there exists a Fr\'{e}chet space  $\mathbb{F}$ containing  $\mathbb{H},$ $\mathbb{F}$ is the completion of $\mathbb{H},$ and for every nonzero vector $f\in \mathbb{H}$ the orbit $\{T^nf,n\geq \}$ meets any open subbase of $\mathbb{F}.$

It was known that the direct sum of two hypercyclic operators need not be hypercyclic but recently De La Rosa and Read~\cite{delaRosa} showed that even the direct sum of a hypercyclic operator with itself $T\oplus T$ need not be hypercyclic.
Bonet and Peris~\cite{BonetPeris} showed that every separable
infinite dimensional Fr\'{e}chet space $\mathbb{F}$ supports a hypercyclic operator.
Moreover, from Shkarin~\cite{Shkarin},
there is a linear operator $T$ such that the direct sum
$T\oplus T\oplus .. .\oplus T = T^{\oplus m}$
of $m$ copies of $T$ is a hypercyclic operator on $\mathbb{F}^m$ for each $m\in \mathbb{N}.$
An $m$-tuple $(T,T,...,T)$ is called {\em disjoint
hypercyclic} if there exists $f\in \mathbb{F}$ such that $(T_1^nf,T_2^nf,...,T_m^nf), n=1,2,...$ is dense in $\mathbb{F}^m.$ See Salas~\cite{Salas} and Bernal-Gonz\'{a}lez~\cite{Bernal} for examples and recent results.

Rezaei~\cite{Rezaei11B} studied weighted composition operators on the space $H(\mathbb{U})$ of holomorphic functions on
$\mathbb{U},$ the open unit disc in $\mathbb{C}.$ Each $\phi\in H(\mathbb{U})$ and holomorphic self-map $\psi$
of $\mathbb{U}$ induce a weighted linear operator $C_{\phi,\psi}$ sending $f(z)$ to $\phi(z)f(\psi(z)).$ This property
 includes both composition $C_\psi,$ ($\phi=1$) and multiplication $M_\phi,$ $\psi=1$) as special cases.
It was shown that any nonzero multiple of $C_\psi$ is chaotic on $H(\mathbb{U})$ if $\psi$ has no fixed point in $\mathbb{U}.$
B\`{e}s et al.~\cite{BesMartinPeris} characterized disjoint hypercyclicity and disjoint supercyclicity of finitely many linear
fractional composition operators (cf. also Bonet and Doma\'{n}ski~\cite{BonetDomanski11} also Zaj\c{a}c~\cite{Zajac12}) acting on spaces of holomorphic functions on the unit
disc, answering a question of Bernal-Gonz\'{a}lez~\cite{Bernal}. Namely, finitely many hypercyclic composition operators $f\mapsto f\circ \varphi$ on the unit disc  $\mathbb{D}$ generated by non-elliptic automorphisms $\varphi$ need not be disjoint nor need they be so on the Hardy space $H^2(\mathbb{D})$ of square-summable power series on the unit disc,
\begin{equation}\label{Hardy}
    H^2(\mathbb{D}) =\left\{ f=z\mapsto \sum_{n=0}^\infty a_nz^n\in H(\mathbb{D}): ||f||^2=\sum_{n=0}^\infty |a_n|^2 <\infty \right\}.
\end{equation}
Shkarin~\cite{Shkarin12} provided an example of a weakly hypercyclic multiplication operator on $H^2(G)$ where $G$ is a region of $\mathbb{C}$ bounded by a smooth Jordan curve $\Gamma$ such that $G$ does not meet the unit ball but $\Gamma$ intersects the unit circle in a non-trivial arc.

Chen and Chu~\cite{ChenChu09,ChenChu11} gave a complete characterization of hypercyclic weighted translation operators on locally compact groups and their homogeneous spaces.
Martin~\cite{Martin11} has notes on hypercyclic properties of groups of linear fractional transformations on the unit disc.
O'Regan and Xian~\cite{OReganX} proved fixed point theorems for maps and multivalued maps between Fr\'{e}chet spaces, using projective limits and the classical Banach theory. Further recent work on set valued maps between Fr\'{e}chet spaces can be found in Galanis et al.{\cite{GBL,GBLP,ORegan} and Bakowska and Gabor~\cite{BakoG}.

Countable products of copies of an infinite-dimensional Banach space are examples of non-normable Fr\'{e}chet spaces that do not admit a continuous norm.
Albanese~\cite{Albanese11} showed that for $\mathbb{F}$  a separable, infinite-dimensional real or complex Fr\'{e}chet
space admitting a continuous norm and $\{v_n\in \mathbb{F}: n \geq 1\}$ a dense set of linearly
independent vectors, there exists a continuous linear operator
$T$ on $\mathbb{F}$ such that the orbit under $T$ of $v_1$ is exactly the set $\{v_n : n \geq 1\}.$  This extended a result of Grivaux~\cite{Grivaux03} for Banach spaces to
the setting of non-normable Fr\'{e}chet spaces that do admit a continuous norm.

\subsection{Semigroups and $n$-tuples of operators}
A Fr\'{e}chet space admits a hypercyclic operator if and only if it
is separable and infinite-dimensional. However, by considering the semigroups generated by multiples of operators, it is possible to obtain hypercyclic behaviour on finite-dimensional spaces.
A semigroup generated by a finite set of $n\times n$ real (or complex) matrices is called {\em hypercyclic} or {\em topologically transitive} if there is a vector with dense orbit in $\mathbb{R}^n$ (or $\mathbb{C}^n$).

Since no finite-dimensional
Banach space admits a hypercyclic operator by Ansari~\cite{Ansari97}, Javaheri~\cite{Javaheri11}  considered a finitely-generated
semigroup of operators instead of a single operator. He gave the following definition as the natural generalization of hypercyclicity to semigroups of operators $\Gamma=\langle T_1,T_2, \ldots, T_k \rangle$ on a finite dimensional vector space over $\mathbb{K}=\mathbb{R} \ {\rm or} \ \mathbb{C}:$ $\Gamma$ is hypercyclic if there exists $x\in\mathbb{K}^n$ such that $\{Tx: T\in \Gamma \}$ is dense in $\mathbb{K}^n.$ Examples were given of $n\times n$ matrices
$A$ and $B$ such that almost every column vector had an orbit under the action of the
semigroup $\langle A, B \rangle$ is dense in $\mathbb{K}^n.$ Costakis et al.~\cite{CostakisHM09}, cf. also~\cite{CostakisHM10}, showed that in every finite dimension there are pairs of commuting matrices which form a locally hypercyclic but non-hypercyclic tuple.
In the non-abelian case, it was shown in \cite{Javaheri11B} that there exists a 2-generator
hypercylic semigroup in any dimension in both real and complex cases. Thus there exists a dense 2-generator semigroup in any dimension in both real and complex
cases. Since powers of a single matrix can never be dense, this result is
optimal.

Ayadi~\cite{Ayadi11} proved that the minimal number of matrices on $\mathbb{C}^n$  required to form
a hypercyclic abelian semigroup on $\mathbb{C}^n$ is $n+1$ and that the
action of any abelian semigroup finitely generated by matrices on $\mathbb{C}^n$
or $\mathbb{R}^n$ is never $k$-transitive for $k \geq 2.$ These answer questions raised by
Feldman and Javaheri~\cite{Javaheri10}.

An {\em $n$-tuple of operators} $T = (T_1, T_2, \ldots , T_n)$ is a finite sequence of length $n$ of commuting continuous linear operators on a locally
convex space $\mathbb{E}$ and $\cal{F} = \cal{F}T$ is the semigroup of strings generated by $T.$ For $f\in \mathbb{E},$ if its orbit under  $\cal{F}$ is dense in $\mathbb{E}$ then the $n$-tuple of operators is called hypercyclic. Feldman~\cite{Feldman08} proved that there are hypercyclic $(n+1)$-tuples of diagonal matrices on $\mathbb{C}^n$
but there are no hypercyclic $n$-tuples of diagonalizable matrices on $\mathbb{C}^n.$
Shkarin~\cite{Shkarin11B} proved that the minimal cardinality of a hypercyclic tuple of
operators is $n+1$ on $\mathbb{C}^n$ and $\frac{n}{2}+\frac{5+(-1)^n}{4}$ on $\mathbb{R}^n$. Also, that there are non-diagonalizable tuples of operators on $\mathbb{R}^2$ which possess an
orbit that is neither dense nor nowhere dense and gave a hypercyclic $6$-tuple of operators on $\mathbb{C}^3$ such that every operator commuting with each
member of the tuple is non-cyclic. A further result was that every infinite-dimensional separable complex (real) Fr\'{e}chet space admits a hypercyclic $6$-tuple ($4$-tuple) $T$ of operators such that there are no cyclic operators commuting with $T.$ Moreover, every hypercyclic tuple $T$ on $\mathbb{C}^2$ or $\mathbb{R}^2$ contains a cyclic operator.

Berm\'{u}dez et al.~\cite{BermudezB&C} investigated hypercyclicity, topological mixing and chaotic maps on Banach spaces. An operator is called mixing if for all nonempty open subsets $U,V$, there is $n\in \mathbb{N}$ such that $T^m(U) \bigcap V \neq \emptyset$ for each $n\geq m.$ An operator is {\em hereditarily hypercyclic} if and only if $T\oplus T$ is hypercyclic~\cite{BesPeris}.  Any hypercyclic operator (on any topological vector space) is transitive. If X is complete separable and metrizable, then the converse implications hold: any transitive operator is hypercyclic and any mixing operator is hereditarily hypercyclic, cf.~\cite{Shkarin11}. Shkarin~\cite{Shkarin11} proved also that a continuous linear operator on a topological vector
space with weak topology is mixing if and only if its dual operator has no finite
dimensional invariant subspaces. Bernal and Grosse-Erdmann~\cite{BernalGE} studied the existence of hypercyclic semigroups of continuous operators on a Banach space. Albanese et al.~\cite{ABR} considered cases when it is possible to  extend Banach space results on $C_0$-semigroups of continuous linear operators to  Fr\'{e}chet spaces. Every operator norm continuous semigroup in a Banach space $X$ has an infinitesimal generator belonging to the space of continuous linear operators on $X;$ an example is given to show that this fails in a general Fr\'{e}chet space. However, it does not fail for countable products of Banach spaces and quotients of such products; these are the Fr\'{e}chet spaces that are quojections, the projective sequence consisting of surjections. Examples include the sequence space $\mathbb{C}^\mathbb{N}$ and the Fr\'{e}chet space of continuous functions $C(X)$ with $X$ a $\sigma$-compact completely regular topological space and compact open topology.

Bayart~\cite{Bayart11} showed that there exist hypercyclic strongly continuous holomorphic groups of operators containing non-hypercyclic operators. Also given were several examples where a family of hypercyclic operators has no {\em common} hypercyclic vector, an important property in linear dynamics, see also Shkarin~\cite{Shkarin10}.

Ayadi et al.~\cite{AyadiMS11} gave a complete characterization of abelian subgroups of $GL(n,\mathbb{R})$ with a locally dense
(resp. dense) orbit in  $\mathbb{R}^n.$ For finitely generated subgroups, this characterization is explicit and it is used
to show that no abelian subgroup of $GL(n,\mathbb{R})$  generated by the integer part of $(n+1/2)$ matrices can have a dense orbit in $\mathbb{R}.$ Several examples are given of abelian groups with dense orbits in $\mathbb{R}^2$ and $\mathbb{R}^4.$ Javaheri~\cite{Javaheri11} gives other results in this context.
Ayadi~\cite{Ayadi11B} characterized hypercyclic abelian affine groups; for finitely
generated such groups, this characterization is explicit. In particular no abelian group generated by $n$ affine maps on $\mathbb{C}^n$ has a dense orbit. An example is given of a group with dense orbit in $\mathbb{C}^2.$

Shkarin~\cite{Shkarin11B} proved that the minimal cardinality of a hypercyclic tuple of
operators on $\mathbb{C}^n$ (respectively, on $\mathbb{R}^n$) is $n+1$ (respectively,
$\frac{n}{2}+\frac{5+(-1)^n}{4}$), that
there are non-diagonalizable tuples of operators on $\mathbb{R}^2$ which possess an
orbit being neither dense nor nowhere dense and construct a hypercyclic $6$-tuple of operators on $\mathbb{C}^3$ such that every operator commuting with each
member of the $6$-tuple is non-cyclic.
It turns out that, unlike for the classical hypercyclicity, there are hypercyclic tuples of operators on finite dimensional spaces. Feldman~\cite{Feldman08}
showed that $\mathbb{C}^n$ admits a hypercyclic $(n+1)$-tuple of operators and
for every tuple of operators on $\mathbb{C}^n$, but not on $\mathbb{R}^n$, every orbit is either dense or is nowhere dense.

The Black-Scholes equation, used (and sometimes misused!~\cite{Stewart12}) for the value of a stock option, yields a semigroup on spaces of continuous
functions on $(0,\infty)$ that are allowed to grow at both $0$ and $\infty,$ which is important since the standard initial value is an unbounded function. Emamirad et al.~\cite{EmamiradGG11}  constructed a family of Banach spaces, parametrized by two market properties on some ranges of which the Black-Scholes semigroup is strongly continuous and chaotic. The proof relied on the Godefroy-Shapiro~\cite{GodefroyShapiro91} Hypercyclicity Criterion, equation (\ref{HC}) above.

\subsection{Topological transitivity and mixing}\label{mixing}
Grosse-Erdmann~\cite{GrE} related hypercyclicity to the topological universality concept, and showed that an operator $T$ is hypercyclic on a separable Fr\'{e}chet space $\mathbb{F}$ if it has the {\em topological transitivity property}: for every pair of nonempty open subsets $U,V\subseteq \mathbb{F}$ there is some $n\in \mathbb{N}$ such that $T^n(U)\bigcap V\neq \emptyset.$ Costakis and V. Vlachou~\cite{CostakisV12} investigated the problem of interpolation by universal, hypercyclic functions. Chen and Shaw~\cite{ChenShaw} linked hypercyclicity to topological mixing, following Costakis and Sambarino~\cite{CS} who showed that if $T^n$ satisfies the Hypercyclicity Criterion then $T$ is {\em topologically mixing} in the sense that:
for every pair of nonempty open subsets $U,V\subseteq \mathbb{F}$ there is some $N\in \mathbb{N}$ such that $T^n(U)\bigcap V\neq \emptyset$ for all $n\geq N.$ Berm\'{u}dez et al.~\cite{BermudezB&C} studied  hypercyclic and chaotic maps on Banach spaces in the context of topological mixing.
See also the summary below on subspace hypercyclicity in \S\ref{subspace} concerning the results of Madore and Mart\'{\i}nez-Avenda\~{n}o~\cite{MadoreMA11} and Le~\cite{Le11}.

B\`{e}s et al.~\cite{BesMartinPeris}  studied mixing and disjoint mixing
behavior of projective limits of endomorphisms of a projective spectrum. In particular, they provided characterization for disjoint hypercyclicity and disjoint supercyclicity of linear
fractional composition operators $C_\varphi:f\mapsto f\circ\varphi$ on $\nu$-weighted Hardy spaces $S_\nu, \ \nu\in \mathbb{R},$ of analytic functions on the unit disc:
\begin{equation}\label{wghtdHardy}
S_\nu =\left\{ f=z\mapsto \sum_{n=0}^\infty a_nz^n\in H(\mathbb{D}): ||f||^2=\sum_{n=0}^\infty |a_n|^2(n+1)^{2\nu} <\infty \right\}
\end{equation}
It was known that a linear fractional composition operator $C_\varphi$ is hypercyclic on $S_\nu$ if and only if $\nu<\frac{1}{2}$ and $C_\varphi$ is hypercyclic on $H^2(\mathbb{D})=S_0,$ equation (\ref{Hardy}). Also, if $\nu<\frac{1}{2}$ then $C_\varphi$ is supercyclic on $S_\nu$ if and only if it is hypercyclic  on $S_\nu.$ B\`{e}s et al.~\cite{BesMartinPeris} extended these results to the projective limit of $\{S_\nu: \nu<\frac{1}{2}\}.$ Zaj\c{a}c~\cite{Zajac12} characterized hypercyclic composition operators in pseudoconvex domains.

Shkarin~\cite{Shkarin11} proved that a continuous linear operator $T$ on a topological vector
space with weak topology is mixing if and only if its dual operator has no finite-dimensional invariant subspace. This result implies the result of Bayart and
Matheron~\cite{BayartMath} that for every hypercyclic operator $T$ on the countable product of copies of $\mathbb{K}=\mathbb{C} \ {\rm or} \ \mathbb{R},$ we have also that $T\oplus T$ is hypercyclic.
Further, Shkarin~\cite{Shkarin11C} described a class of topological vector spaces admitting a mixing uniformly continuous
operator group $\{T_t\}_{t\in \mathbb{C}^n}$ with holomorphic dependence on the parameter $t,$ and a class of topological vector spaces admitting no supercyclic
strongly continuous operator semigroups $\{T_t\}_{t\geq0}.$

\subsection{Subspace hypercyclicity}\label{subspace}
Madore and Mart\'{\i}nez-Avenda\~{n}o~\cite{MadoreMA11} introduced the concept
of {\em subspace hypercyclicity}: a continuous linear operator $T$ on a Hilbert space $\mathbb{H}$ is {\em $M$-hypercyclic} for a subspace
$M$ of $\mathbb{H}$ if there exists a vector such that the intersection of its orbit and $M$
is dense in $M.$ Those authors proved several results analogous to the hypercyclic case.
For example, if $T$ is subspace-hypercyclic, then its
spectrum must intersect the unit circle, but not every element of the spectrum need do so; subspace-hypercyclicity is a strictly infinite-dimensional
phenomenon; neither compact operators nor hyponormal (i.e. $||Tx||\geq ||T^*x||, \forall x\in \mathbb{H}$) bounded operators are subspace-hypercyclic.

For closed $M$ in separable Banach $\mathbb{E},$ Madore and Mart\'{\i}nez-Avenda\~{n}o~\cite{MadoreMA11} showed that $M$-hypercyclicity is implied by $M$-transitivity---i.e. for all disjoint nonempty open subsets $U,V$ of $M$ there is a number $n$ such that $U\bigcap T^{-n}V$ contains a nonempty open set of $M.$
Le~\cite{Le11}  gave a sufficient condition
for $M$-hypercyclicity and used it to show that it need not imply $M$-transitivity.

Desch and Schappacher~\cite{DeschSchappacher11} defined the (weakly) topological transitivity of a semigroup $\cal{S}$  of bounded linear operators on a real Banach space as the property
for all nonempty (weakly) open sets $U,V$ that for some $T\in \cal{S}$ we have $TU\bigcap V \neq \emptyset.$
They characterized
weak topological transitivity of the families of operators $\{S^t | t\in\mathbb{N}\},\ \{kS^t | t\in\mathbb{N}, k>0\},$
and $\{kS^t | t\in\mathbb{N}, k\in\mathbb{R}\},$ in terms of the point spectrum of the dual operator $S^*$ cf. also~\cite{BakkaliT12}.
Unlike topological transitivity in the norm topology, which is equivalent to
hypercyclicity with concomitant highly irregular behaviour of the semigroup, Desch and Schappacher~\cite{DeschSchappacher11} illustrated quite good behaviour of weakly topologically transitive semigroups. They gave an example using the positive-definite bounded self-adjoint
$$S:L^2([0,1])\rightarrow L^2([0,1] : u(\xi) \mapsto \frac{u(\xi)}{\xi+2}.$$
Then $S=S^*$ and has empty point spectrum so $\{S^t : t\in \mathbb{N} \}$ is weakly topologically transitive but cannot be weakly hypercyclic because $S^t\rightarrow 0$ in
operator norm if $t\rightarrow\infty.$ They point out that weak transitivity is in fact a weak property. For, a weakly open set in an infinite-dimensional Banach
space contains a subspace of finite codimension but an apparently small neighborhood
contains many large vectors, easily hit by trajectories.

Rion's thesis~\cite{Rion11} is concerned particularly with hypercyclicity of the Aluthge transform of weighted shifts on $\ell^2(Z).$  In Chapter 4 he considered also the distribution of hypercyclic vectors over the
range of a hypercyclic operator, pointing out that if $x$ is a hypercyclic vector for $T,$ then so is $T^nx$ for all $n\in \mathbb{N},$ and $T^nx$ is in the range of $T.$
Since moreover, the range of $T$ is dense, one might expect that most if not all of an operators hypercyclic vectors lie in its range. However, Rion~\cite{Rion11} showed for every non-surjective hypercyclic operator $T$ on a Banach
space, the set of hypercyclic vectors for $T$ that are not in its range is large, in that it is not expressible as a countable union of nowhere dense sets, providing also a sense by which the range of an arbitrary hypercyclic operator $T$ is large in its set of
hypercyclic vectors for $T.$

\subsection{Chaotic behaviour}
A continuous linear operator $T$ on a topological vector space $\mathbb{E}$ has a {\em periodic point} $f\in \mathbb{E}$ if, for some $n\in \mathbb{N}$ we have $T^nf=f.$ The operator
$T$ is {\em cyclic} if for some $f\in \mathbb{E}$ the span of $\{T^nf, n\geq 0 \}$ is dense in $\mathbb{E}.$ On finite-dimensional spaces there are many cyclic operators but no hypercyclic operators.
The operator $T$ is called {\em chaotic}~\cite{GrE11} if it is hypercyclic and its set of periodic points is dense in $\mathbb{E}.$ Each operator on the Fr\'{e}chet space of analytic functions on $\mathbb{C}^N,$ which commutes with all translations and is not a scalar multiple of the identity, is chaotic~\cite{GodefroyShapiro91}.

Rezaei~\cite{Rezaei11B} investigated weighted composition operators on the space $H(\mathbb{U})$ of holomorphic functions on
$\mathbb{U},$ the open unit disc in $\mathbb{C}.$ Each $\phi\in H(\mathbb{U})$ and holomorphic self-map $\psi$
of $\mathbb{U}$ induce a weighted linear operator $C_{\phi,\psi}$ sending $f(z)$ to $\phi(z)f(\psi(z)).$
It was shown that any nonzero multiple of $C_\psi$ is chaotic on $H(\mathbb{U})$ if $\psi$ has no fixed point in $\mathbb{U}.$
Berm\'{u}dez et al.~\cite{BermudezB&C} studied  hypercyclic and chaotic maps on Banach spaces in the context of topological mixing.
Emamirad et al.~\cite{EmamiradGG11}  constructed a family of Banach spaces, parametrized by two market properties on
some ranges of which the Black-Scholes semigroup is strongly continuous and chaotic. That proof relied on the
Godefroy-Shapiro~\cite{GodefroyShapiro91} Hypercyclicity Criterion, equation (\ref{HC}) above.

The conjugate set ${\cal{C}}(T)=\{L^{-1} T L: L \ {\rm invertible}\}$ of a hypercyclic operator
$T$ consists entirely of hypercyclic operators, and those hypercyclic
operators are dense in the algebra of bounded linear operators with respect to the strong operator topology.
Chan and Saunders~\cite{ChanSanders} showed that, on an infinite-dimensional Hilbert space,  there is a path of chaotic operators, which is dense in the operator algebra with the strong operator topology, and along which every operator has the exact same dense $G_\delta$ set of hypercyclic vectors.
Previously~\cite{ChanSanders09} they showed that the conjugate set of any
hypercyclic operator on a separable, infinite dimensional Banach space always contains
a path of operators which is dense with the strong operator topology, and yet the set
of common hypercyclic vectors for the entire path is a dense $G_\delta$ set. As a corollary, the hypercyclic
operators on such a Banach space form a connected subset of the operator algebra with the strong operator topology.

 Originally defined on a metric space $(X,d),$ a Li-Yorke chaotic map $f:X\rightarrow X$ is such that there exists an uncountable subset $\Gamma\subset X$ in which every pair of distinct
 points $x,y$ satisfies
 \begin{equation}\label{chaotic}
   \liminf_n d(f^nx,f^ny) =0 \ \ {\rm and} \ \limsup_n d(f^nx,f^ny) > 0
 \end{equation}
 then $\Gamma$ is called a {\em scrambled} set. The map $f$ is called {\em distributionally chaotic} if there is an $\epsilon>0$ and an uncountable set $\Gamma_\epsilon\subset X$ in which every pair of distinct points $x,y$ satisfies
  \begin{equation}\label{chaotic}
   \liminf_{n\rightarrow\infty}\frac{1}{n}|\{k: d(f^kx,f^ky)<\epsilon, 0\leq k<n\}| =0  \ {\rm and} \ \limsup_{n\rightarrow\infty}\frac{1}{n}|\{k: d(f^kx,f^ky)<\epsilon, 0\leq k<n\}| =1
 \end{equation}
 and then $\Gamma_\epsilon$ is called a distributionally $\epsilon$-scrambled set, cf. Mart\'{\i}nez-Gim\'{e}nez~\cite{MartinezOP09}.
For example, every hypercyclic operator $T$
on a Fr\'{e}chet space $F$ is Li-Yorke chaotic with respect to any (continuous) translation invariant metric: just fix a hypercyclic vector $x$ and $\Gamma = \{\lambda x: |\lambda|\leq 1\}$ is a scrambled set for $T.$
Berm\'{u}dez et al.~\cite{BermudezBMP11} characterized on Banach spaces continuous Li-Yorke chaotic bounded linear operators $T$ in terms of the existence of irregular vectors; here, $x$, is irregular for $T$ if
 \begin{equation}\label{irreg}
   \liminf_n ||T^nx|| =0 \ \ {\rm and} \ \limsup_n ||T^nx|| = \infty.
 \end{equation}
Sufficient `computable' criteria for Li-Yorke chaos were given, and they established some additional conditions for the existence of dense scrambled sets. Further, every infinite dimensional separable Banach space was shown to admit a distributionally chaotic operator which is also hypercyclic, but from Mart\'{\i}nez-Gim\'{e}nez et al.~\cite{MartinezOP09}there are examples of backward shifts on K\"{o}the spaces of infinite-dimensional matrices which are uniformly distributionally chaotic and not hypercyclic. K\"{o}the spaces provide a natural class of Fr\'{e}chet sequence spaces (cf. also Golinsky~\cite{Golinski12}) in which many typical examples of weighted shifts are chaotic.
Mart\'{\i}nez-Gim\'{e}nez et al.~\cite{MartinezOP12} showed that neither hypercyclicity nor the 
mixing property is a sufficient condition for distributional chaos. 

 The existence
of an uncountable scrambled set in the Banach space setting may not be as strong an indication of complicated dynamics as in the compact metric space case. For example, it may happen that the span of a single vector becomes an uncountable
scrambled set~\cite{BermudezBMP11}. This led Subrahmonian Moothathu~\cite{SubMoothathu11} to look for some feature stronger than uncountability for a scrambled set in the
Banach space setting. He showed that if an operator is hypercyclic, so it admits a vector with dense orbit, then it has a scrambled set in the strong sense of requiring
linear independence of the vectors in the scrambled set.

Following the Chen and Chu~\cite{ChenChu09,ChenChu11} complete characterization of hypercyclic weighted translation operators on locally compact groups and their homogeneous spaces, Chen~\cite{CCChen11} then characterized chaotic weighted translations, showing that the density of periodic points of a weighted translation implies hypercyclicity. However, a weighted translation operator is not hypercyclic if it is generated by a group element of finite order~\cite{ChenChu11}. A translation operator is never chaotic because its norm cannot exceed unity, but a weighted translation can be chaotic.
It was known that for a unimodular complex number $\alpha$ the rotation $\alpha T$ of a hypercyclic operator on a complex Banach space is also hypercyclic but Bayart and Bermudez~\cite{BayartBermudez09} showed that on a separable Hilbert space there is a chaotic operator $T$ with $\alpha T$ not chaotic. Chen~\cite{CCChen11} proved that this is not the case for chaotic weighted translation operators becuase their rotations also are chaotic.

Desch et al.~\cite{DeschSW} gave a sufficient condition  for a stronly continuous semigroup of
bounded linear operators on a Banach space to be chaotic in terms of the spectral properties
of its infinitesimal generator, and studied applications to several differential equations with constant coefficients.
 Astengo and Di Blasio~\cite{AstengoBlasio11} extended this study to the chaotic and hypercyclic behaviour of the strongly continuous  modified heat semigroup of operators generated by perturbations of the Jacobi Laplacian with a multiple of the identity on $L^p$ spaces.

\end{document}